\documentclass[12pt]{article}
\usepackage{amssymb,amsmath}
\newtheorem{thm}{Theorem}
\newtheorem{cor}[thm]{Corollary}
\newtheorem{lemma}[thm]{Lemma}
\newenvironment{defin}{\medskip\noindent{\sc
Definition}.}{\goodbreak\medskip}
\newenvironment{nota}{\medskip\noindent{\sc
Notations}.}{\goodbreak\medskip}
\newenvironment{remk}{\noindent{\sc
Remark}.}{\goodbreak\vskip10pt}
\newtheorem{prop}[thm]{Proposition}

\def\demo{\medskip\goodbreak\noindent
     \hbox{\sc Proof \kern .3em}\ignorespaces}%
  \def \qedbox{$\square$}%
  \def \qed{\hglue1mm\hfill{\ifmmode\qedbox
     \else\unskip\ \hglue0mm\hfill\qedbox\medskip
      \goodbreak\fi}}%
\def\enddemo{\qed\goodbreak\vskip10pt}%

\def\qed{\hglue1mm\hfill\raise -2pt\hbox{\vrule\vbox to 10pt{\hrule width
4pt
                  \vfill\hrule}\vrule}}

\newcommand{\T}{\mathbb {T}}

\newcommand{\esse}{\mathbb {S}}
\newcommand{\R}{\mathbb {R}}

\newcommand{\N}{\mathbb {N}}
\newcommand{\Nc}{\mathcal {N}}

\newcommand{\Vc}{\mathcal {V}}

\newcommand{\Cc}{\mathcal {C}}
\newcommand{\Ic}{\mathcal {I}}

\newcommand{\Mc}{\mathcal {M}}
\newcommand{\Kc}{\mathcal {K}}
\newcommand{\Ec}{\mathcal {E}}
\newcommand{\Gc}{\mathcal {G}}
\newcommand{\Lc}{\mathcal {L}}

\newcommand{\Sc}{\mathcal{S}}

\renewcommand{\Im}{\mathrm{Im}}

\begin{document}
\title{{ Green bundles, Lyapunov exponents and regularity along  the supports of the minimizing measures  }}
\author{M.-C. ARNAUD
\thanks{ANR Project BLANC07-3\_187245, Hamilton-Jacobi and Weak KAM Theory}
\thanks{ANR DynNonHyp}
\thanks{Universit\'e d'Avignon et des Pays de Vaucluse, Laboratoire d'Analyse non lin\' eaire et G\' eom\' etrie (EA 2151),  F-84 018Avignon,
France. e-mail: Marie-Claude.Arnaud@univ-avignon.fr}
}
\maketitle
\abstract{   In this article, we study the minimizing measures of the Tonelli Hamiltonians. More precisely, we study the relationships between the so-called Green bundles and various notions as:
\begin{enumerate}
\item[$\bullet$] the Lyapunov exponents of minimizing measures;
\item [$\bullet$]   the weak KAM solutions.
\end{enumerate}
In particular, we deduce that the support of every minimizing measure $\mu$, all of whose Lyapunov exponents are zero, is $C^1$-regular $\mu$-almost everywhere.
}
\medskip


\begin{center}
{\bf R\'esum\'e}
\end{center}
Dans cet article, on \'etudie les mesures minimisantes de Hamitoniens de Tonelli. Plus pr\'ecis\'ement, on explique quelles relations existent entre les fibr\'es de Green et diff\'erentes notions comme~:
\begin{enumerate}
\item[$\bullet$] les exposants de Lyapunov des mesures minimisantes;
\item[$\bullet$] les solutions KAM faibles.
\end{enumerate}
On en d\'eduit par exemple que si tous les exposants de Lyapunov d'une mesure mi\-ni\-misante $\mu$ sont nuls, alors le support de cette mesure est $C^1$-r\'egulier en $\mu$-presque tout point.
\medskip

\noindent Mots clefs: Orbites et mesures minimisantes, exposants de Lyapunov, th\'eorie KAM faible, fibr\'es de Green, r\'egularit\'e des solutions de l'\'equation de Hamilton-Jacobi.

\medskip

\newpage
\tableofcontents
\newpage
\section{Introduction}
In this article, $M$ is a closed  $n$-dimensional manifold and $\pi~: T^*M\rightarrow M$ its cotangent bundle. We consider a Tonelli Hamiltonian $H~: T^*M\rightarrow \R$,  i.e.  a $C^2$ function that is  strictly $C^2$-convex  and superlinear in the fiber. The Hamiltonian flow associated with such a function is denoted by $(\varphi_t)_{t\in\R}$ or $(\varphi_t^H)_{t\in\R}$. To such a Hamiltonian, there corresponds a Lagrangian function $L~: TM\rightarrow \R$ that has the same regularity as $H$ and is also superlinear and strictly convex in the fiber. The corresponding Euler-Lagrange flow is denoted by $(f_t)_{t\in\R}$.

\medskip

For such a Hamiltonian system, it is usual to study its ``minimizing objects''; more precisely, a   piece of orbit   $(\varphi_t(q,p))_{t\in [a,b]}=(q_t, p_t)_{t\in[a,b]}$ is minimizing if the arc $(q_t)_{t\in[a, b]}$  minimizes the action functional $A_L$ defined by $A_L(\gamma )=\int_a^bL(\gamma (t), \dot\gamma (t))dt$ among the $C^2$-arcs joining $q_a$ to $q_b$. More generally, if $I$ is an interval and $(\varphi_t)_{t\in I}=(q_t, p_t)_{t\in I}$ is an orbit piece, we say that it is minimizing if for every segment $[a,b]\subset I$, its restriction to $[a,b]$ is minimizing. Then we call the set of  points of $T^*M$ whose (complete) orbit is minimizing the {\em Ma\~n\'e set}. We denote it by $\Nc^*(H)$ and its projection, the {\em projected Ma\~n\'e set}, is denoted by: $\Nc (H)=\pi (\Nc^*(H))$. The Ma\~n\'e set is non empty, compact and invariant by the Hamiltonian flow (see \cite{Fa1}). The first proof of the non-emptiness of the Ma\~n\'e set is due to J.~Mather: he proved in the 90's in \cite{mather1} the existence of minimizing measures.

We are interested in invariant subsets of the Ma\~n\'e set, i.e.   subsets that are the union of some minimizing orbits. More precisely, we would like to know if we can say something about the regularity of such subsets (we will be more precise very soon. It's a kind of differentiability) and particularly if there is a link between the dynamic of the flow restricted to such a  set and the regularity of the set.

\medskip

The oldest result in this direction concerns the time-dependent case~: considering a symplectic twist map of the annulus $T^*\esse$, G.~Birkhoff proved in the 1920's that any essential invariant curve is the graph of a Lipschitz map (see \cite{Bir1} or \cite{He1}). It is easy to prove that such a curve is action minimizing. In the case of   higher dimensions, M.~Herman proved in \cite{He2} that any $C^0$-Lagrangian graph of $T^*\T^n$  that is invariant  by a symplectic  twist map is, in fact,  the graph of a Lipschitz map. A related result in the autonomous case is that any $C^1$-Hamilton-Jacobi solution of a Tonelli Hamiltonian is, in fact, $C^{1,1}$ (see \cite{Fa2}). As Rademacher's theorem says to us that any Lipschitz function is differentiable Lebesgue almost everywhere,  these results are a kind of regularity result.

\medskip

In \cite{Arna2}, we did, in fact, improve these results of regularity in the autonomous case, proving that  if a $C^0$-Lagrangian graph is invariant by a Tonelli flow, and if one of the two following hypotheses is satisfied:
\begin{enumerate}
\item[$\bullet$] $\dim M=2$ and all the singularities of $H$ are non degenerate;
\item[$\bullet$] the dynamic of the restriction of the flow to the invariant graph is Lipschitz conjugate to a translations' flow;
\end{enumerate}
 then the invariant graph is, in fact, $C^1$ almost everywhere (this is stronger than just differentiable). 
Let us point out  that any of the two previous hypotheses implies that the dynamic of the restricted flow to the graph is soft on a certain sense (our arguments are not very precise, but we only want to give a certain intuition of the forthcoming result); indeed, when $\dim M=2$, if we reduce the dynamic modulo the vector field, we obtain a 1-dimension dynamic, and it is known at least in the differentiable case that the Lyapunov exponents of a dynamic on the circle are zero. The same is true for any dynamic that is Lipschitz conjugate to a translation.  

We gave a similar results for the invariant curves of the twist maps of the annulus in \cite{Arna1}, proving that Birkhoff's result can be improved: any essential invariant curve of a symplectic twist map of the annulus $T^*\esse$ is the graph of a Lipschitz map that is $C^1$ Lebesgue almost everywhere.

\medskip

Hence, it seems reasonable to try to find a relationship between the Lyapunov exponents of any  minimizing measure and the regularity of its support, where an invariant  measure is {\em minimizing} if its support is in the Ma\~n\'e set. 

For a twist map of the annulus $T^*\esse$,  we studied the ergodic minimizing measures in \cite{Arna3} and proved that the $C^1$-regularity (we will be more precise very soon) of its support is equivalent to the fact that the Lyapunov exponents are zero. Hence, in a certain way, in this case,  ``$C^1$-irregularity'' is equivalent to non-vanishing Lyapunov exponents.

\medskip

The question that we ask now ourselves is the following: what  can  we say for higher dimensions? Is the irregularity (in a sense we will soon specify) of the support of a minimizing ergodic measure equivalent to non-vanishing exponents?

\medskip

A first and obvious answer is: no. Indeed, let us consider the following example: $(\psi_t)$ is an Anosov   flow defined on the cotangent bundle $T^*\Sc$ of a closed surface $\Sc$. Let $\Nc=T^*_1\Sc$ be its unitary cotangent bundle, which is a 3-manifold invariant by $(\psi_t)$. Then a method due to Ma\~n\'e (see \cite{Man1}) allows us to define a Tonelli Hamiltonian $H$ on $T^*\Nc$ such that the restriction of its flow $(\varphi_t)$ to the zero section $\Nc$ is $(\psi_t)$: the Lagrangian $L$ associated with $H$  is defined by: $L(q,v)=\frac{1}{2}\| \dot\psi (q)-v\|^2$ where $\| .\|$ is any Riemannian metric on $\Nc$. In this case, the zero section is very regular (even $C^\infty$), but the Lyapunov exponents of every invariant measure whose support is contained in $\Nc$ are non zero (except two, the one corresponding to the flow direction and the one corresponding to the energy direction). Hence, it may happen that some exponents are non zero and the support of the measure is very regular\dots\\

In fact, the other implication is true: we will see that the nullity of the Lyapunov exponents implies the regularity of the support of the considered measure.\\
Let us now explain in a detailed way in which kind of regularity we are  interested: 

\begin{defin}
Let $A$ be a subset of a manifold $M$ and let $a$ belong to $A$. The contingent cone to $A$ at $a$ is the set of the tangent vectors $v\in T_aM$ such that there exist   a sequence $(a_n)$ of elements of $A$ and   a sequence $(t_n)$ of positive real numbers  such that (we write everything in a chart, but this is independent of the chosen chart):
$$\lim_{n\rightarrow \infty} \frac{1}{t_n}(a_n-a)=v.$$
We denote it by: $\Cc_aA$.
\end{defin}
This notion of contingent cone is due to Bouligand (see \cite{Bou}). The contingent cone is never empty (it always contains the null vector), and it is equal to the null vector if, and only if, $a$ is an isolated point of $A$.

We   will see later that the sets in which we are interested are contained in some (weak) Lagrangian manifolds. Our definitions of $1$-regularity and $C^1$-regularity seems very natural for such sets:  

\begin{defin}
Let $A$ be a subset of a symplectic manifold $M$ and let $a$ belong to $A$. We say that $A$ is $1$-regular at $a$ if the contingent cone to $A$ at $a$ is contained in a Lagrangian subspace of $T_aM$. \\
We say that $A$ is $C^1$ regular at $a$ if there exists a Lagrangian subspace $\Lc$ of $T_aM$ such that: for every sequence $(a_n, v_n\in \Cc_{a_n}A)$ such that $\displaystyle{\lim_{n\rightarrow \infty}a_n=a}$ and the sequence $(v_n)$ converges to an element $v$ of $T_aM$, then $v\in \Lc$.
\end{defin}

Let us notice that this notion of $C^1$-regularity is slightly different from  the ones given in \cite{Arna1}, \cite{Arna2} and \cite{Arna3}: the notions given in these former articles are a little stronger. This notion of $C^1$-regularity  is stronger than the notion of $1$-regularity, which   is nothing else but the notion of differentiability for the $C^0$- Lagrangian graphs (see \cite{Arna2} for a definition of $C^0$-Lagrangian graphs).  \\

The measures that we   study are the minimizing ones, that is the ones that are invariant and whose supports are contained in the Ma\~{n}\'e set. Then we prove:

\begin{thm}\label{th1}
Let $H~: T^*M\rightarrow \R$ be a Tonelli Hamiltonian and let $\mu$ be an ergodic minimizing probability measure all of  whose   Lyapunov are zero. Then, at $\mu$-almost every point of the support ${\rm supp}(\mu)$ of $\mu$, the set ${\rm supp}(\mu)$ is $C^1$-regular.
\end{thm}

Hence:
\begin{enumerate}
\item[$\bullet$] we succeed in proving that a kind of ``soft dynamic'' implies some $C^1$-regularity;
\item[$\bullet$] we know that we can have simultaneously a strong dynamic (for example hyperbolic)  and a $C^\infty$-regularity.
\end{enumerate}

In fact, we obtain more precise results than this theorem; for example, an interesting question is: what happens if there are simultaneously some zero and non-zero exponents?\\
To explain what happens, we need to introduce some other notions. Let us begin by  recalling what the Green bundles are. These Lagrangian bundles were introduced by L.~Green  in 1958 in \cite{Gr} for geodesic flows to prove some rigidity results.  For the existence and the construction of these bundles, the reader is referred to \cite{Arna2}, \cite{C-I} or \cite{It1}. We recall:

\begin{defin} Here, $V(x)=\ker D\pi (x)$ designates the linear vertical.\\
Let $(\varphi_t(q,p))_{t\in ]-\infty, 0]}$ be a minimizing negative orbit; then the positive  Green bundle $G_+$ is defined along this orbit by: $\displaystyle{G_+(x)=\lim_{t\rightarrow +\infty} D\varphi_t.V(\varphi_{-t}x)}$.\\
Let $(\varphi_t(q,p))_{t\in [0, +\infty[}$ be a minimizing positive orbit; then the negative  Green bundle $G_-$ is defined along this orbit by: $\displaystyle{G_-(x)=\lim_{t\rightarrow +\infty} D\varphi_{-t}.V(\varphi_{t}x)}$.
\end{defin}
Hence, at every point of the Ma\~n\'e set, the two Green bundles are defined.

Let us recall that the two Green bundles are Lagrangian, invariant under the linearized flow $D\varphi_t$, transverse to the vertical, that they depend semi-continuously on the considered point (see \cite{Arna2} for the definition of semi-continuity of Lagrangian subspaces transverse to the vertical), that $G_-\leq G_+$ (see \cite{Arna2} for the definition of the order between two planes transverse to the vertical; in coordinates, this corresponds to the usual order on the set of symmetric matrices whose  Lagrangian  subspaces are the graphs.). Hence, if $\mu$ is an ergodic  minimizing probability measure,    the integer $\dim(G_-(x)\cap G_+(x))$ is constant $\mu$ almost everywhere.

We obtain a result linking the dimension of the intersection of the two Green bundles to the number of non zero Lyapunov exponents:

\begin{thm}\label{thdimlyap}
Let $H~: T^*M\rightarrow \R$ be a Tonelli Hamiltonian and let $\mu$ be an ergodic minimizing probability measure. Then the two following assertions are equivalent:
\begin{enumerate}
\item[$\bullet$] at $\mu$ almost every point, $\dim (G_-(x)\cap G_+(x))=p$;
\item[$\bullet$] $\mu$ has exactly $2p$ zero Lyapunov exponents, $n-p$ positive ones and $n-p$ negative ones.
\end{enumerate}
\end{thm}

Let us mention some former related results: 
\begin{enumerate}
\item[$\bullet$] in \cite{C-I}, the authors prove that the transversality of the two Green bundles along an energy level implies that the restriction of the flow to this level is Anosov; they use some ideas about quasi-Anosov dynamics due to R.~Ma\~n\'e that are contained in \cite{Man2}; in \cite{Eb1}, P.~Eberlein gives the same statement for the geodesic flows;
\item[$\bullet$]  we proved in \cite{Arna3} that  any quasi-hyperbolic symplectic cocycle above a compact  set is hyperbolic; we can apply this result to any minimizing compact invariant subset $K$ contained in an energy level $\Ec$ without singularity: considering the restricted/reduced dynamical system to the energy level $\Ec$ modulo the vector-field (see \cite{Arna2} p 899 for the construction), we deduce that the transversality of the Green bundles in the energy level above $K$ is equivalent to the partial hyperbolicity of the linearized flow along $K$ with a center bundle's dimension equal to 2;

\item[$\bullet$] concerning the non-uniform case (i.e. the case of minimizing measures), the only known result was a formula giving the entropy due to A.~Freire \& R.~Ma\~n\'e (see \cite{F-M}). Roughly speaking, by integrating some functional along one of the two Green bundles, they compute the sum of the positive Lyapunov exponents. This formula was generalized in \cite{C-I} to any Tonelli Hamiltonian. But this formula doesn't say to us how many non-zero Lyapunov exponents exist: it  only gives   the sum of the positive Lyapunov exponents. Let us mention too that G.~Knieper gives a nicer formula in his (non-published) thesis.
\end{enumerate}
 
 To prove theorem  \ref{th1}, we  recall in section \ref{sec3} some points of the recent weak KAM theory developped by A.~Fathi in \cite{Fa1}. In this section too, we give some statements concerning the relationships between weak KAM solutions and the Green bundles. We don't give them in the introduction because we would need all the notions that will be defined in section \ref{sec3}, but the interested reader can go to section \ref{sec3}. Roughly speaking, the theorem asserts that along the support of the minimizing measures, the contingent cones to the weak KAM pseudographs is not far from some cone delimited by the two Green bundles. \\
 
 Theorem \ref{thdimlyap} is proved in section \ref{sec2}. The statement concerning the relationships between the weak KAM solutions and the Green bundles are contained in section \ref{sec3} and the proofs are in section \ref{sec4}.

 \section{Green bundles and Lyapunov exponents}\label{sec2}
 
 In this section, we   prove theorem \ref{thdimlyap}.
We consider an ergodic minimizing measure $\mu$ that is not  the Dirac measure at a critical point and we denote the integer such that we have $\mu$ almost everywhere: $\dim G_-\cap G_+=p$ by $p$. Let us recall the dynamical criterion  that is proved in \cite{Arna2}:

\begin{prop}[dynamical criterion] \label{dyncrit} Let $(x_t)$ be a minimizing and relatively compact orbit. Let $v\in T_{x_0}(T^*M)$. Then:\\
-- if $v\notin G_-(x_0)$, then $\displaystyle{\lim_{t\rightarrow +\infty}\| D\pi\circ D\varphi_t.v\|=+\infty}$;\\
--if $v\notin G_+(x_0)$, then $\displaystyle{\lim_{t\rightarrow +\infty}\| D\pi\circ D\varphi_{-t}.v\|=+\infty}$.
\end{prop}
and some direct consequences of this criterion: 
 
\begin{remk} 1) We deduce from the dynamical criterion that  the Hamiltonian vector-field  $X_H$ belongs to the two Green bundles. This implies that $p\geq 1$. Because these two Green bundles are Lagrangian, this implies that $G_+$ and $G_-$ are tangent to the Hamiltonian levels $\{ H=c\}$. \\
2) Moreover, we deduce also that if there is an Oseledet splitting (this will be   precisely defined very soon)  $T(T^*M)=E^s\oplus E^c\oplus E^u$ above a minimizing compact set $K$, then $E^s\subset G_-$ and $E^u\subset G_+$. Because the flow is symplectic, $E^u$ and $E^s$ are isotropic and orthogonal to $E^c$ for the symplectic form (see \cite{bochi-viana1}). Moreover, $E^{s\bot}=E^s\oplus E^c$ (where $\bot$ designates the orthogonal subspace for the symplectic form) and $E^{u\bot}=E^u\oplus E^c$; we deduce that: $G_-(x)=G_-(x)^\bot\subset  E^{s\bot}=E^s\oplus E^c$ and similarly that $G_+(x)\subset E^u(x)\oplus E^c(x)$. Hence, finally:
 $$E^s(x)\subset G^-(x)\subset E^s(x)\oplus E^c(x)\ {\rm and}\ E^u(x)\subset G^+(x)\subset E^u(x)\oplus E^c(x)$$
 and then: $G_-(x)\cap G_+(x)\subset E^c(x)$. Hence, $G_-\cap G_+$ being an isotropic subspace of the symplectic subspace $E^c$, we obtain: $\dim E^c\geq 2\dim (G_-\cap G_+)$. The dimension of the intersection of the two Green bundles gives a lower bound to the number of zero Lyapunov  exponents. Theorem \ref{thdimlyap} says to us that this inequality is, in fact, an equality. 
Let us notice that when $p=n$, we directly have the conclusion of the theorem because $\dim E^c\geq 2\dim M$ implies that $\dim E^c=2n$.\\
We have the same results for a hyperbolic or partially hyperbolic dynamic. Let us notice that in the hyperbolic case, $G_-$ (resp. $G_+$) is nothing else but the stable (resp. unstable) bundle $E^s$ (resp. $E^u$)\\
3) Let us consider the case of a K.A.M. torus that is a graph (when $M=\T^n$):  the dynamic on this torus is $C^1$ conjugated to a flow of irrational translations on the torus $\T^n$; M.~Herman proved in \cite{He2} that such a torus is Lagrangian, and it is well-known that any invariant Lagrangian graph is locally minimizing. Then the orbit of every vector tangent to the K.A.M. torus is bounded, and  belongs to $G_-\cap G_+$. In this case, the two Green bundles are equal to the tangent space to the invariant torus.
\end{remk}
Let us introduce some notations:

 \begin{nota} Oseledet's theorem implies that there exist an invariant subset $N$ of $T^*M$ with full $\mu$-measure, some real numbers  $0<\lambda _1  <\lambda_2 <\dots < \lambda_{q }$ and a (measurable) splitting with constant dimensions above $N$: 
 $$T_x(T^*M)=E^s_1(x)\oplus E^s_2(x)\oplus\dots\oplus E^s_{q }(x)\oplus E^c(x)\oplus E^u_1(x)\oplus E^u_2(x)\oplus \dots \oplus E^u_{q }(x)$$ such that:
 \begin{enumerate}
 \item[$\bullet$] for every $v\in E^s_j(x)\backslash\{ 0\}$; $\displaystyle{\lim_{t\rightarrow\pm \infty}\frac{1}{t}\log\left( \| D\varphi_t(x)v\| \right)=-\lambda_j}$;
 \item[$\bullet$] for every $v\in E^c(x)\backslash\{ 0\}$; $\displaystyle{\lim_{t\rightarrow\pm \infty}\frac{1}{t}\log\left( \| D\varphi_t(x)v\| \right)=0}$;
\item[$\bullet$] for every $v\in E^u_j(x)\backslash\{ 0\}$; $\displaystyle{\lim_{t\rightarrow\pm \infty}\frac{1}{t}\log\left( \| D\varphi_t(x)v\| \right)=+\lambda_j}$.\\
We may ask, too,  that: $\forall x\in N, \dim(G_-(x)\cap G_+(x))=p$.
 \end{enumerate}
\end{nota}

Let us recall that the stable bundle $E^s(x)=E^s_1(x)\oplus E^s_2(x)\oplus\dots\oplus E^s_{q }(x)$ and the unstable one $E^u(x)=E^u_1(x)\oplus E^u_2(x)\oplus\dots\oplus E^u_{q }(x)$ are isotropic (for the symplectic form) and that $E^c(x)$ is a symplectic subspace of $T_x(T^*M)$ that is orthogonal (for $\omega$) to $E^s(x)\oplus E^u(x)$.  Moreover, we have: $\dim E^s_i=\dim E^u_i$.

\subsection{Reduction of the problem}
As in the statement of   theorem  \ref{thdimlyap},  we assume that  $\mu$ is a minimizing ergodic measure whose support is not reduced to a point and  that $p\in [1, n]$ is so that at $\mu$-almost every point $x$, the intersection of  the Green bundles $G_+(x)$ and $G_-(x)$ is $p$-dimensional. We deduce from the previous remark that for every $x\in N$:  $G_+(x)\cap G_-(x)\subset E^c(x)$ and $E^s(x)\oplus E^u(x)= \left( E^c(x)\right)^\bot\subset G_+(x)^\bot+ G_-(x)^\bot=G_-(x)+ G_+(x)$.  

\begin{nota}
We introduce the two notations:  $E(x)=G_-(x)+G_+(x)$ and    $R(x)= G_-(x)\cap G_+(x)$. We denote the reduced space: $F(x)=E(x)/R(x)$  by $F(x)$ and we denote   the canonical projection $p~: E\rightarrow F$ by $p$.  As $G_-$ and $G_+$ are invariant by the linearized flow $D\varphi_t$, we may define a reduced cocycle $M_t~: F\rightarrow F$. But $(M_t)$ is not continuous, because $G_-$ and $G_+$ don't vary continuously.\\
Moreover, we introduce the notation: $\Vc(x)=V(x)\cap E(x)$ is the trace of the linearized vertical on $E(x)$ and $v(x)=p(\Vc(x))$ is the projection of $\Vc (x)$ on $F(x)$. We introduce a notation for the images of the reduced vertical $v(x)$ by $M_t$: $g_t(\varphi_tx)=M_tv(x)$.
\end{nota}

The subspace $E(x)$ of $T_x(T^*M)$ is  co-isotropic with $E(x)^{\bot}=R(x)$. Hence $F(x)$ is nothing else than the symplectic space that is obtained by symplectic reduction of $E(x)$. We denote its symplectic form by $\Omega$. Hence we have: $\forall (v, w)\in E(x)^2, \Omega (p(v), p(w))=\omega (v,w)$. Moreover,  $(M_t)$ is a symplectic cocycle.\\
We can notice, too, that $\dim E(x)=\dim (G_-(x)+G_+(x))=\dim G_-(x)+\dim G_+(x)-\dim (G_-(x)\cap G_+(x))=2n-p$ and deduce that $\dim F(x)=\dim E(x)-\dim (G_-(x)\cap G_+(x))=2(n-p)$.

\begin{nota} If $L$ is any Lagrangian subspace of $T_x(T^*M)$, we denote $(L\cap E(x))+R(x)$ by $\tilde{L}$ and  $p(\tilde{L})$ by $l$. 
\end{nota}

\begin{lemma}\label{bernard1} If $L\subset T_x(T^*M)$ is Lagrangian, then $\tilde{L}$ is also Lagrangian and $l=p(\tilde L)=p(L\cap E(x))$ is a Lagrangian subspace of $F(x)$. Moreover, $p^{-1}(l)=\tilde{L}$ . In particular, $v(x)$ is a Lagrangian subspace of $F(x)$ and $p^{-1}(v(x))=\Vc(x)+R(x)$.

\end{lemma}
\demo
We just have to prove that $\tilde L$ is Lagrangian, the other assertions being easy consequences of this fact. \\
We begin by proving that $\tilde L$ is isotropic. If $u, u'\in L\cap E(x)$ and $v, v'\in R(x)$, then $\omega (u+v, u'+v')=0$ because $L$ is Lagrangian and then $\omega (u, u')=0$ and because $R(x)\subset E(x)^\bot$.\\
Let us determine $\dim \tilde{L}$. Let $L'$ be such that: $L=(E(x)\cap L)\oplus L'$. Then the dimension of  $L\cap R(x)=(L+E(x))^\bot$ is: $2n-(\dim L+E(x))=2n-(2n-p+\dim L')=p-\dim L'$. We deduce: $\dim \tilde L=\dim (L\cap E(x))+\dim R(x)-\dim (L\cap R(x))=\dim (L\cap E(x))+p-(p-\dim L')=\dim (L\cap E(x))+\dim L'=\dim L$.
\enddemo

\begin{lemma}\label{L2}
The subspace $v(x)$ is a Lagrangian subspace of $F(x)$. Moreover, for every $t\not=0$, $g_t(\varphi_tx)=M_tv(x)$ is transverse to $v(\varphi_t(x))$
\end{lemma}
\demo
The first sentence is contained in lemma \ref{bernard1}.\\
  Let us  consider $t\not=0$ and let us assume that $ M_tv(x) \cap v(\varphi_tx)\not=\{ 0\}$.  We may assume that $t>0$ (or we replace $x$ by $\varphi_t(x)$ and $t$ by $-t$).
  
  Then there exists $v\in \Vc(x)\backslash \{  0\}$ such that$D\varphi_t(x)v\in \Vc (\varphi_t x)+(G_-(\varphi_tx)\cap G_+(\varphi_tx))$.  Let us write $D\varphi_t(x)v=w+g$ with $w\in \Vc (\varphi_tx)$ and $g\in R(\varphi_tx)$.  We know that the orbit has no conjugate vector (because the measure is minimizing); hence $g\not=0$. 
  
  Moreover, we proved in  \cite{Arna2} that $D\varphi_t V(x)$ is strictly above $G_-(\varphi_tx)$, i.e. that:
  $$\forall h\in G_-(\varphi_tx), \forall k\in V(\varphi_tx),    h+k\in D\varphi_t V(x)\backslash \{ 0\}\Rightarrow  \omega  (h, h+k)> 0.$$
  We deduce that: $\omega (g, w+g)>0$. 
  
 This contradicts: $D\varphi_t(x)v\in E(\varphi_tx)=\left(G_+(\varphi_tx)\cap G_-(\varphi_tx)\right)^\bot\subset (\R g)^\bot$.
\enddemo
As in \cite{Arna2}, we ask ourselves   what the order between the different Lagrangian subspaces $g_t(x)=M_tv(\varphi_{-t}x)$ is. Let us recall how we define this order:

\begin{defin}
Let $g_1$ and $g_2$ be two subspaces of $F(x)$ that are transverse to the (reduced) vertical $v(x)$. Let $f(x)=F(x)/v(x)$ be the reduced space and $P(x)~: F(x)\rightarrow f(x)$ the canonical projection. Then to every $w\in f(x)$, we can associate a unique $\ell_1(w)\in g_1$ (resp. $\ell_2(w)\in g_2$) such that: $P( \ell_1(w))= w$ (resp. $P(\ell_2(w))= w$). We then define the altitude of $g_2$ above $g_1$, which is a quadratic form defined on $f(x)$, by: $q(g_1,g_2)(w)=\Omega (\ell_1(w), \ell_2(w))$.\\
We say that $g_2$ is above (resp. strictly above) $g_1$ when $q(g_1, g_2)$ is positive semi-definite (resp. positive definite). We write $g_1\leq g_2$ (resp. $g_1<g_2$).
\end{defin}

\begin{lemma}\label{bernard2} Let $L_1$, $L_2$ be two Lagrangian subspaces of $T_x(T^*M)$ transverse  to $V(x)$ such that at least one of them is contained in $E(x)$. Then, if $L_1<L_2$ (resp. $L_1\leq L_2$), we have: $l_1$ and $l_2$ are transverse to $v(x)$ and $l_1 < l_2$ (resp. $l_1\leq l_2$). We deduce that $p(G_-)<p(G_+)$.

\end{lemma}
\demo
We assume that $L_2\subset E(x)$ and that $L_1<L_2$. Let $v_1\in L_1\cap E(x)$ be a non-zero vector of $L_1\cap E(x)$. As $L_1$ and $L_2$ are transverse to $V(x)$, there exists a unique $v_2\in L_2$ such that $v_2-v_1\in V(x)$. Moreover, as $v_1, v_2\in E(x)$, we have $v_2-v_1\in \Vc (x)$  and $p(v_2)-p(v_1)\in v(x)$. Hence:
$$\Omega (p(v_1), p(v_2))=\omega (v_1, v_2)>0.$$
This means exactly that $l_1<l_2$.\\
To deduce the assertion for $\leq$, we can use a limit.\\
As $G_-\leq G_+$, we deduce that $p(G_-)\leq p(G_+)$. Because of the definition of $E(x)$, $R(x)$ and $F(x)$, $p(G_-)$ and $p(G_+)$ are transverse and  then  $p(G_-)<p(G_+)$.

\enddemo
 
 \begin{lemma}
 If $\mu$ is a minimizing measure, for every $x\in{\rm supp}\mu$, for all $0<t<s$, we have:
 $$g_{-t}(x)<g_{-s}(x)<g_s(x)<g_t(x).$$
 \end{lemma}
  \demo
  The map $(t\in \R^*\rightarrow g_t(x))$ is continuous; moreover, we know by lemma \ref{L2} that if $t\not=s$, then $g_t(x)$ is transverse to $g_s(x)$. Hence, the index of $q(g_s(x), g_t(x))$ is constant for $(s, t)\in\Ec$ where $\Ec$ is one of the sets: $\{ (s,t); 0<s<t\}$; $\{ (s,t); s<0<t\}$, $\{ (s,t); s<t<0\}$. Hence, we only have to determine this index for one point $(s,t)$ of each of these three sets.\\
 We  prove the result only for the first set, the other inequalities being very similar. Let us fix $s>0$ and introduce the notation $G_s(x)=D\varphi_sV(\varphi_{-s}x)$. Then  $\tilde G_s(x)$ is a Lagrangian subspace of $E(x)$ that is transverse to the vertical because $\tilde G_s(x)\cap V(x)=\tilde G_s(x)\cap \Vc (x)=(\tilde G_s(x)\cap \tilde V(x))\cap \Vc(x)=p^{-1}(g_s(x)\cap v(x))\cap \Vc(x)=R(x)\cap \Vc(x)=\{ 0\}$. We assume that $t>0$ is very small and we work in a chart, with   symplectic coordinates defined in \cite {Arna2} (p 897)  such that  the ``horizontal'' subspace of $T_x(T^*M)$ is  $G_-(x)$. A vector of $ G_t(x)=D\varphi_t(\varphi_{-t}x)V(\varphi_{-t}x)$ is $(h, S^+_t(x)h)$ and it is proved in \cite{Arna2} (p 894) that $S^+_{t}(x)\sim \frac{1}{t}D$ where $D$ is a fixed positive definite matrix. Hence, for $t>0$ small enough, we have $\tilde G_s<G_t$. We deduce from lemma \ref{bernard2} that $g_s=p(\tilde G_s)<p(G_t)=g_t$.\\
  \enddemo
 \begin{defin}
 As in \cite{Arna2}, when $t$ tends to $\pm  \infty$, we find two $M_t$-invariant Lagrangian sub-bundle of $F(x)$ that are: $\displaystyle{g_-(x)=\lim_{t\rightarrow -\infty} g_t(x)}$ and $\displaystyle{g_+(x)=\lim_{t\rightarrow +\infty} g_t(x)}$; they are transverse to $v(x)$ and satisfy: $g_-(x)\leq g_+(x)$. We call them the reduced Green bundles.
\end{defin}
\remk Then we have: $\forall t>0, g_{-t}(x)<g_-(x)\leq g_+(x)<g_t(x)$.  If we use the notations $\tilde G_\pm(x)=p^{-1}(g_\pm(x))$, then $\tilde G_\pm$ are transverse to the vertical because $\tilde G_\pm(x) \cap V(x)=\tilde G_\pm (x)\cap \Vc(x)=(\tilde G_\pm\cap \tilde V(x))\cap \Vc (x)=p^{-1}(g_\pm (x)\cap v(x))\cap \Vc(x)=R(x)\cap \Vc(x)=\{ 0\}$.
Moreover, $\tilde G_-(x)\leq \tilde G_+(x)$ and the two bundles $\tilde G_-$, $\tilde G_+$, are invariant by the linearized flow $(D\varphi_t)$. Theorem 3.11 of \cite{Arna2} asserts that any invariant Lagrangian bundle that is transverse to the vertical is between the two Green bundles. We deduce that  $G_-(x)\leq \tilde G_-(x)\leq \tilde G_+(x)\leq G_+(x)$. We can then use lemma \ref{bernard2} and we obtain: $p(G_-(x))\leq g_-(x)\leq g_+(x)\leq p(G_+(x))$.
\begin{lemma}\label{derlem}
We have:  $\forall x\in{\rm supp}\mu, g_-(x)=p(G_-(x))<p(G_+(x))= g_+(x) $.
\end{lemma}
\demo
Because of the last remark, we just have to prove that on ${\rm supp} \mu$: $  g_- \leq p(G_- )<p(G_+ )\leq g_+  $. Because of lemma \ref{bernard2}, we just have to prove that   $  g_- \leq p(G_- )$ and $p(G_+ )\leq g_+  $. But $p(G_\pm)$ is a lagrangian subspace of $F(x)$ whose orbit is transverse to the vertical. We can use a similar statement to proposition 3.11 of \cite{Arna2} to deduce the inequalities.
\enddemo
Hence we have proved that $\tilde G_\pm=G_\pm$, the notation $\tilde G_\pm$ will disappear from tnow on.

 \subsection{Reduced Green bundles and   Lyapunov exponents}

 We have to  be careful because the bundles that we consider are not continuous and, as this is noted in \cite{Arna2},  we don't use a continuous change of coordinates, but just  a bounded one when we say that $G_-$ or $G_+$ is the horizontal subspace (the matrix $P$ that is necessary to change the coordinates is uniformly bounded, as $P^{-1}$). \\
 We choose at every point $x\in N$ some (linear) symplectic coordinates $(Q,P)$ of $F(x)$ such that $v(x)$ has for equation: $Q=0$ and $g_+(x)$ has for equation $P=0$. We will be more precise on this choice later. Then the matrix of $M_t(x)$ in these coordinates is a symplectic matrix: $M_t(x)=\begin{pmatrix}
 a_t(x)&b_t(x)\\
0&d_t(x)\\
 \end{pmatrix}$. As $M_t(x)v(x)=g_t(\varphi_tx)$ is a Lagrangian subspace of $E(\varphi_tx)$ that is transverse to the vertical,  then $\det b_t(x)\not=0$ and there exists a symmetric matrix $s_t^+(\varphi_tx)$ whose graph is $g_t(\varphi_tx)$, i.e: $d_t(x)=s_t^+(\varphi_t(x))b_t(x)$. Moreover, the family $(s^+_t(x))_{t>0}$ being decreasing and tending to zero (because by hypothesis the horizontal is $g_+$), the symmetric matrix $s^+_t(\varphi_tx)$ is positive definite. Moreover, the matrix $M_t(x)$ being symplectic, we have: 
 $$\left(M_t(x)\right)^{-1}=\begin{pmatrix} {}^td_t(x)&-{}^tb_t(x)\\
 0&{}^ta_t(x)\\
\end{pmatrix}$$
and by definition of $g_{-t}(x)$, if it is the graph of the matrix $s^-_t(x)$ (that is negative definite), then: ${}^ta_t(x)=-s_t^-(x){}^tb_t(x)$ and finally:
$$M_t(x)=\begin{pmatrix}
-b_t(x)s_t^-(x)& b_t(x)\\
0& s_t^+(\varphi_tx)b_t(x)\\
\end{pmatrix}
$$
Let us be now more precise in the way we choose our coordinates; we may associate  an almost complex structure $J$ and then a Riemannian metric $(.,.)_x$ defined by: $(v,u)_x=\omega (x)(v, Ju)$ with the symplectic form $\omega$ of $T^*M$; from now on,  we work with this fixed Riemannian metric of $T^*M$. We choose on $ G_+(x)=p^{-1}(g_+(x))$ an orthonormal basis whose last vectors are in $R(x)$ and complete it in a symplectic base whose last vectors are in $V(x)$. We denote  the associated  coordinates of $T_x(T^*M)$ by $(q_1, \dots, q_n, p_1, \dots, p_n)$. These (linear) coordinates don't depend in a continuous way on the point $x$ (because $G_+$ doesn't), but in a bounded way. Then $G_-(x)=p^{-1}(g_-(x))$ is the graph of a symmetric matrix whose kernel is $R(x)$ and then on $G_-(x),$ we have: $p_{n-p+1}=\dots =p_n=0$. An element of $E(x)$ has coordinates such that $p_{n-p+1}=\dots =p_n=0$, and an element of $F(x)=E(x)/R(x)$ may be identified with an element with coordinates $(q_1, \dots , q_{n-p}, 0, \dots , 0, p_1, \dots , p_{n-p}, 0, \dots  , 0)$. We then use on $F(x)$ the norm $\displaystyle{\sum_{i=1}^{n-p}(q_i^2+p_i^2)}$, which is the norm for the Riemannian metric   of the considered element of $F(x)$. Then this norm depends in a measurable way on $x$. 

Let us now notice the following fact: $\mu$ being ergodic for the flow  $(\varphi_t)$, there exists a dense $G_\delta$ subset $A$ of $\R$ such that, for every $t\in A$, the diffeomorphism $\varphi_t$ is ergodic. As it is   simpler for us to work with a diffeomorphism  instead of a flow, we fix such a $t\in A$.  We assume that $t=1$ (if not we replace $H$ by $\frac{1}{t}H$).\\

\begin{lemma}\label{LJ}
For every $\varepsilon >0$, there exists a measurable subset $J_\varepsilon$ of $N$ such that:
\begin{enumerate}
\item[$\bullet$] $\mu (J_\varepsilon)\geq 1-\varepsilon$;
\item[$\bullet$] on $J_\varepsilon$, $(s_n^+)$ and $(s_n^-)$ converge uniformly ;
\item[$\bullet$] there exists two constants $\beta=\beta(\varepsilon)>\alpha=\alpha(\varepsilon)>0$ such that: $\forall x\in J_\varepsilon,\beta{\bf 1}\geq  -s_-(x)\geq \alpha {\bf 1}$ where $g_-$ is the graph of $s_-$.
\end{enumerate}
\end{lemma}
\demo  This is a consequence of Egorov theorem and of the fact that on $N$, $g_+$ and $g_-$ are transverse and then $-s_-$ is positive definite.
\enddemo
 We deduce:
 \begin{lemma}\label{LCVU}
 Let $J_\varepsilon$ be as in the previous lemma. On the set $\{ (n,x)\in\N\times J_\varepsilon, \varphi_n(x)\in J_\varepsilon\}$, the sequence of conorms $(m(b_n(x))$ converge uniformly to $+\infty$, where $m(b_n)=\| b_n^{-1}\| ^{-1}$.
 \end{lemma}
 \demo Let $n, x$ be as in the lemma.\\
 The matrix  $M_n(x)=\begin{pmatrix}
 -b_n(x)s_n^-(x)& b_n(x)\\
 0& s_n^+(\varphi_nx)b_n(x)\\
 \end{pmatrix}
 $ being symplectic, we have: \\
 $-s_n^-(x){}^tb_n(x)s_n^+(\varphi_nx)b_n(x)={\bf 1}$ and thus  
 $-b_n(x)s_n^-(x){}^tb_n(x)s_n^+(\varphi_nx)={\bf 1}$ and:\\
  $b_n(x)s_n^-(x){}^tb_n(x)=-\left(s_n^+(\varphi_nx)\right)^{-1}$. \\
 We know that on $J_\varepsilon$, $(s_n^+)$ converges uniformly to zero. Hence,  for every $\delta>0$, there exists $N=N(\delta) $ such that: $n\geq N\Rightarrow \| s_n^+(\varphi_nx)\|\leq \delta$. Moreover, we know that $\| s_n^-(x)\|\leq \beta$. Hence, if we choose $\delta'=\frac{\delta^2}{\beta}$, for every $n\geq N=N(\delta')$ and $x\in J_\varepsilon$ such that $\varphi_nx\in J_\varepsilon$, we obtain: 
 $$\forall v\in\R^p,\beta \| {}^tb_n(x)v\|^2= {}^tv b_n(x)(\beta{\bf 1}){}^tb_n(x)v\geq - {}^tv b_n(x)s_n^-(x){}^tb_n(x)v={}^tv\left(s_n^+(\varphi_nx)\right)^{-1}v$$
 and we have: ${}^tv\left(s_n^+(\varphi_nx)\right)^{-1}v\geq \frac{\beta}{\delta^2}\| v\|^2$ because $s_n^+(\varphi_nx)$ is a positive definite matrix that is less than $\frac{\delta^2}{\beta}{\bf 1}$. We finally obtain: $\| {}^tb_n(x)v\|\geq \frac{1}{\delta}\| v\|$ and then the result that we wanted.
 \enddemo
 From now we fix  a small constant $\varepsilon>0$, associate  a set $J_\varepsilon$ with $\varepsilon$  via lemma \ref{LJ} and two constants $0<\alpha<\beta$; then  there exists $N\geq 0$ such that
 $$\forall x\in J_\varepsilon, \forall n\geq N, \varphi_n(x)\in J_\varepsilon\Rightarrow m(b_n(x))\geq \frac{2}{\alpha}.$$

 \begin{lemma}
 Let $J_\varepsilon$ be as in lemma \ref{LJ}. For $\mu$-almost point $x$ in $J_\varepsilon$, there exists a sequence of integers $(j_n)=(j_n(x))$ tending to $+\infty$ such that: 
 $$\forall n\in \N, m(b_{j_n}(x)s_{j_n}(x))\geq \left( 2^\frac{1-\varepsilon}{2N}\right)^{j_n}.$$
 \end{lemma} 
 \demo
 As $\mu$ is ergodic for $\varphi_1$, we deduce from Birkhoff ergodic theorem that for almost every point $x\in J_\varepsilon$, we have:
 $$\lim_{\ell\rightarrow +\infty}\frac{1}{\ell}\sharp \{ 0 \leq k\leq \ell-1; \varphi_k(x)\in J_\varepsilon\}=\mu (J_\varepsilon)\geq 1-\varepsilon.$$
 We introduce the notation: $N(\ell)=\sharp \{ 0 \leq k\leq \ell-1; \varphi_k(x)\in J_\varepsilon\}$.\\
 For such an $x$ and every $\ell\in\N$, we find a number $n(\ell)$ of integers:
 $$0=k_1\leq k_1+N\leq k_2\leq k_2+N\leq k_3\leq k_3+N\leq \dots \leq k_{n(\ell)}\leq \ell$$
 such that $\varphi_{k_i}(x)\in J_\varepsilon$ and $n(\ell)\geq [\frac{N(\ell)}{N}]\geq \frac{N(\ell)}{N}-1$. In particular, we have: $\frac{n(\ell)}{\ell}\geq\frac{1}{N}(\frac{N(\ell)}{\ell}-\frac{N}{\ell})$, the right term converging to $\frac{\mu (J_\varepsilon)}{N}\geq \frac{1-\varepsilon}{N}$ when $\ell$ tends to $+\infty$. Hence, for $\ell$ large enough, we find: $n(\ell)\geq 1+  \ell \frac{1-\varepsilon}{2N}$.\\ 
As $\varphi_{k_i}(x)\in J_\varepsilon$ and $k_{i+1}-k_i\geq N$, we have:  $m(b_{k_{i+1}-k_i}(\varphi_{k_i}(x)))\geq \frac{2}{\alpha}$. Moreover, we have: $m(s_{k_{i+1}-k_i}^-(\varphi_{k_i}x))\geq \alpha$; hence: 
$$m(b_{k_{i+1}-k_i}(\varphi_{k_i}x)s_{k_{i+1}-k_i}^-(\varphi_{k_i}x))\geq 2.$$
But the matrix $-b_{k_{n(\ell)}}(x)s^-_{k(n(\ell))}(x)$ is the product of $n(\ell)-1$
 such matrix. Hence:
 $$m(b_{k_{n(\ell)}}(x)s^-_{k(n(\ell))}(x))\geq 2^{n(\ell)-1}\geq 2^{\ell\frac{1-\varepsilon}{2N}}\geq \left( 2^\frac{1-\varepsilon}{2N}\right)^{k_{n(\ell)}}.
 $$
 
 \enddemo
Let us now come back to the whole tangent space $T_x(T^*M)$ with a slight  change in the coordinates that we use. We defined the symplectic coordinates $(q_1, \dots , q_n, p_1, \dots , q_n)$ and now we use the non symplectic ones: \\
$(Q_1, \dots, Q_n,P_1, \dots , P_n)=(q_{n-p+1}, \dots, q_n, q_1, \dots, q_{n-p}, p_1, \dots , p_n)$.  Then:
\begin{enumerate}
 \item[$\bullet$] $(Q_1, \dots , Q_p)$ are coordinates in $R(x)$;
 \item[$\bullet$] $(Q_1, \dots , Q_n)$ are coordinates in $G_+(x)$;
 \item[$\bullet$] $(Q_1, \dots , Q_n, P_{1}, \dots , P_{n-p})$ are coordinates of $E(x)=G_+(x)+G_-(x)$.
 \end{enumerate}
We write then the matrix of $D\varphi_t(x)$  in these  coordinates $(Q_1, \dots , Q_n, P_1, \dots  s , P_n)$ (which are not symplectic):
$$\begin{pmatrix}
A^1_t(x)&A^2_t(x)&A^3_t(x)&A^4_t(x)\\
0&b_t(x)s_t^-(x)&b_t(x)&A^5_t(x)\\
0&0& s_t^+(\varphi_tx)b_t(x)&A^6_t(x)\\
0&0&0&A^9_t(x)\\
\end{pmatrix}$$
where the blocks correspond  to the decomposition $T_x(T^*M)=E_1(x)\oplus E_2(x)\oplus E_3(x)\oplus E_4(x)$ with $\dim E_1(x)=\dim E_4(x)=p$ and $\dim E_2(x)=\dim E_3(x)=n-p$.\\
We have noticed that $E_1(x)=E(x)\subset E^c(x)$ and that $G_+(x)=E_1(x)\oplus E_2(x)$.\\
If $x\in J_\varepsilon$, we have found a sequence $(j_n)$ of integers tending to $+\infty$ so that: 
 $$\forall n\in \N, m(b_{j_n}(x)s^-_{j_n}(x))\geq \left( 2^\frac{1-\varepsilon}{2N}\right)^{j_n}.$$
 We deduce:
 $$\forall v\in E_2(x)\backslash \{ 0\}, \frac{1}{j_n}\log\left( \| b_{j_n}(x)s^-_{j_n}(x)v\|\right)\geq  \frac{1-\varepsilon}{2N}\log 2 +\frac{\|v\|}{j_n};$$
 and because $E_1(x)\subset E^c(x)$: 
 $$\forall v\in G_+(x)\backslash E_1(x), \liminf_{n\rightarrow \infty}\frac{1}{n}\log \| D\varphi_n(x)v\|\geq \frac{1-\varepsilon}{2N}\log 2.$$
Hence there are at least $n-p$  Lyapunov exponents  bigger than $ \frac{1-\varepsilon}{2N}\log 2$ and then bigger than $0$ for the linearized flow. Because this flow is symplectic, we deduce that it has at least $n-p$ negative Lyapunov exponents (see \cite{bochi-viana1}). As we noticed that the linearized flow has at least $2p$ zero Lyapunov exponents, we deduce that $\mu$ has   exactly  $n-p$ positive Lyapunov  exponents, exactly  $n-p$ negative Lyapunov exponents and  exactly $2p$ zero Lyapunov exponents.\\
This finishes the proof of theorem \ref{thdimlyap}.

\remk   Let us notice that we proved too that for $x\in N$ (i.e. generic in the Oseledet's sense), we have:  $E^u(x)\subset G_+(x)$, and then $G_+(x)=E^u(x)\oplus R(x)  $
 \section{Weak K.A.M. solutions and Green bundles}\label{sec3}
 In this section, we   recall the weak KAM theory and give a relationship between some tangent cones to the pseudographs of the weak KAM solutions and the Green bundles. These results imply theorem \ref{th1}. The proofs are given in section \ref{sec4}.
 
 \subsection{Weak KAM theory} We don't give any proof in this section, but all the results that we give are proved in \cite{Fa1} or \cite{Be1}.
 
 \begin{nota}
  If $t>0$, the function $A_t~: M\times M\rightarrow \R$ is defined by: 
$$A_t(q_0, q_1)=\inf_\gamma \int_0^tL(\gamma (s), \dot\gamma (s))ds=\min_\gamma \int_0^tL(\gamma (s), \dot\gamma (s))ds$$
where the infimum is taken on the set of $C^2$ curves  $\gamma~: [0, t]\rightarrow M$ such that $\gamma (0)=q_0$ and $\gamma (t)=q_1$. 
 \end{nota}
\begin{defin}
\begin{enumerate}
\item A function $v~: V\rightarrow \R$ defined on a subset $V$ of $\R^d$ is {\em $K$-semi-concave} if for every $x\in V$, there exists a linear form $p_x$ defined on $\R^d$ so that:
$$\forall y\in V, v(y)\leq v(x)+p_x(y-x)+K\| y-x\|^2.$$
Then we say that $p_x$ is a {\em K-super-differential} of $v$ at $x$.
\item Let us fix a finite atlas ${\cal A}$ of the manifold $M$; a function $u~: M\rightarrow \R$ is {\em $K$-semi-concave} if for every chart $(U, \phi)$ belonging to ${\cal A}$, $u\circ\phi^{-1}$ is $K$-semi-concave. Then a {\em $K$-super-differential} of $u$ at $q$ is $p_x\circ D\phi(q)$ where $p_x$ is a $K$-super-differential of $u\circ \phi^{-1}$ at $x=\phi(q)$.
\end{enumerate}
A semi-concave function is always Lipschitz and then differentiable almost  everywhere  and for such a function, we define its pseudograph: 
a {\em pseudograph} is the graph ${\cal G}(du)$ of $du$, where $u~: M\rightarrow \R$ is a semi-concave function. \\
A function $u~: M\rightarrow \R$ is $K$-semi-convex if $-u$ is $K$-semi-concave. We have a notion of sub-differential and the anti-pseudograph of a semi-convex function $u$ is $\Gc (du)$.
\end{defin}
 It is proved in \cite{Be1} that $A_t$ is semi-concave and that for every minimizing curve $\gamma~: [0, t]\rightarrow M$ between $q_0$ and $q_1$,  $(-\frac{\partial L}{\partial v}(\gamma (0), \dot\gamma (0)), \frac{\partial L}{\partial v}(\gamma (t), \dot\gamma (t)))$ is a super-differential of $A_t$ at $(q_0, q_1)$. It is proved, too, that $A_t(.,q_1)$ is differentiable at $q_0$ if, and only if, $A_t(q_0,.)$ is differentiable at $q_1$ if, and only if, there exists a unique minimizing curve $\gamma~: [0, t]\rightarrow M$ joining $q_0$ to $q_1$.\\

 We denote the two Lax-Oleinik semi-groups associated with $L$  by $(T_t)_{t>0}$ and $(\breve T_t)_{t>0}$; for $u\in C^0(M, \R)$ , they are  defined by:
 $$T_tu(q)=\min_{q'\in M} (u(q')+A_t(q',q))\ {\rm and}\ \breve T_tu(q)=\max_{q'\in M} (u(q')-A_{t}(q, q'))$$
 A function $u~: M\rightarrow \R$ is a negative (resp. positive) weak KAM solution if there exists $c\in\R$ such that: $\forall t>0, T_tu=u-ct$ (resp. $\forall t>0, \breve T_t u=u+ct$).
 
Then there   exist at least one positive and one negative   weak K.A.M. solutions  (see \cite{Fa1} or \cite{Be1}).  The constant $c$ is   unique and is called Ma\~n\'e's critical value. If $u_-$ is a negative weak KAM solution  and $u_+$ a positive one, then $u_-$ is semi-concave and $u_+$ is semi-convex. 
Let us introduce the Mather set:

\begin{defin}
The Mather set, denoted by $\Mc^*(H)$,  is the union of the supports of the minimizing measures. The projected Mather set is $\Mc(H)=\pi(\Mc^*(H))$.
\end{defin}
J.~Mather proved that $\Mc^*(H)$ is compact, non-empty and that it is a Lipschitz graph above a compact part of the zero-section of $T^*M$.

A.~Fathi proved in \cite{Fa1} that if $u_-$ is a negative weak KAM solution,  there exists a unique positive weak KAM solution $u_+$ such that $u_{-|\Mc(H)}=u_{+|\Mc(H)}$. Such a pair $(u_-, u_+)$ is called a pair of conjugate weak KAM solutions. For such a pair, we have:
  \begin{enumerate}
  \item[$\bullet$] $\forall q\in \Mc(H), u_-(q)=u_+(q)$; let us denote the set of equality: $\Ic (u_-, u_+)=\{ q; u_-(q)=u_+(q)\}$  by $\Ic (u_-,u_+)$; then $\Mc(H)\subset \Ic(u_-, u_+)$; 
  \item[$\bullet$] $u_-$ and $u_+$ are differentiable at every point $q\in \Ic(u_-, u_+)$; for such a $q$ we have $(q, du_-(q))\in \Nc^*(H)$; when $q\in\Mc (H)$ and  $(q, p)\in  \Mc^*(H)$ is its lift to $\Mc^*(H)$, then $  du_-(q)=  du_+(q)=p$;
  \item[$\bullet$] $u_+\leq u_-$.
  \end{enumerate}
  Moreover, it is proved in \cite{Be1} that if $q$ is a point of differentiability of $T_tu$ (resp. $\breve T_tu$), then the minimum (resp. maximum) in the definition of $T_tu(q)$ (resp. $\breve T_tu$) is attained at a unique $q'$ and there is a unique curve $\gamma~: [0, t]\rightarrow M$  minimizing between $q'$ and $q$ (resp. $q$ and $q'$);  in this case: $\frac{\partial L}{\partial v}(q, \dot\gamma (t))=dT_tu(q)$ (resp. $\frac{\partial L}{\partial v}(q, \dot\gamma (0))=d\breve T_tu(q)$).

 \subsection{Comparison between the weak KAM solutions and the Green bundles }
 If $(u_-, u_+)$ is a pair of conjugate weak KAM solutions, if $q\in \Ic (u_-, u_+)$, we have seen that $(q, du_-(q))=(q, du_+(q))\in\Nc^*(H)$. Hence, the two Green subspaces $G_-(q, du_-(q))$ and $G_+(q, du_+(q))$ exist. Let us introduce two other Lagrangian subspaces:
 
 \begin{nota}
 If the orbit of $x$ is minimizing, if $G_-(x)$  is the graph of the symmetric matrix $s_-(x)$ and   $G_+(x)$ the graph of the symmetric matrix  $s_+(x)$, we denote the graph of $\tilde s_-(x)= 2s_-(x)-s_+(x)$ (resp. $\tilde s_+(x)= 2s_+(x)-s_-(x)$) by $\tilde G_-(x)$ (resp. $\tilde G_+(x)$).\\
 If $\Delta s(x)=  s_+(x)- s_-(x)$, then $\Delta s (x)$ is positive semi-definite and we have: $\tilde s_-=s_- -\Delta s$ and $\tilde s_+=s_++\Delta s$.\\
 Moreover, if $s$ is a positive semi-definite matrix, we will denote by $p_s$ the orthogonal projection on its image $\Im (s)$ and by $\Lambda (s)$ is greatest eigenvalue: $\Lambda (s)=\| s\|$.
 \end{nota}
 Let us notice that   $G_-(x)=G_+(x)$ if, and only if, $\tilde G_-(x)=G_-(x)=G_+(x)=\tilde G_+(x)$. Moreover, we always have: $\tilde G_-(x)\leq G_-(x)\leq G_+(x)\leq \tilde G_+(x)$. The bundle $\tilde G_-$ is lower semi-continuous and the bundle $\tilde G_+
$ is upper semi-continuous, and they are continuous at the points where $G_-=G_+$.

Let us recall that if $x\in A\subset T^*M$, $\Cc_xA$ designates the contingent cone to $A$ at $x$, that was defined in the introduction. 
 \begin{thm}\label{greenkam}
 Let $(u_-, u_+)$ be a pair of conjugate weak KAM solutions and let $q$ belong to $\Ic(u_-, u_+)$. Then we have:  $\forall (X, Y)\in  \Cc_{(q, du_-(q))}\Gc (du_-),  $
 $$\| Y-\tilde s_-(q, du_-(q))X\|\leq 2\sqrt{\|\Delta s(q, du_-(q))\|} .\sqrt{\Delta s(q, du_-(q))(X,X)}$$
 $$\leq 2\Lambda (\Delta s(q, du_-(q))). \| p_{\Delta s(q, du_-(q))}(X)\|$$
 and:   $\forall (X, Y)\in  \Cc_{(q, du_+(q))}\Gc (du_+),  $
 $$\| Y-\tilde s_+(q, du_+(q))X\|\leq 2\sqrt{\|\Delta s(q, du_+(q))\|} .\sqrt{\Delta s(q, du_+(q))(X,X)}$$
 $$\leq 2\Lambda (\Delta s(q, du_+(q))). \| p_{\Delta s(q, du_+(q))}(X)\|$$
 
 \end{thm}
 We postpone the proof of this theorem to section \ref{sec4}.\\
 As $\Mc^*(H)\subset \Gc(du_-)\cap \Gc(du_+)$, we deduce: 
 \begin{cor}\label{cormes}
 If $x$ is an element of $\Mc^*(H)$, then we have:  $\forall (X, Y)\in \Cc_x\Mc^*(H)$, 
 $$ \max \{ \| Y-\tilde s_-(x)X\|, \| Y-\tilde s_+(x)X\|\}\leq 
 2\sqrt{\|\Delta s(x)\|} .\sqrt{\Delta s (x) (X,X)}\leq 2\Lambda (\Delta s(x)). \| p_{\Delta s(x)}(X)\|$$
 \end{cor}
 Now, we use theorem \ref{thdimlyap}: if $\mu$ is an ergodic minimizing measure whose Lyapunov exponents are zero, then we have $\mu$-almost everywhere: $G_-=G_+$ i.e. $\Delta s=0$. We deduce from corollary \ref{cormes} that $\Cc_x({\rm supp} \mu) \subset G_-(x)=G_+(x)$ at $\mu$ almost every point. This implies that ${\rm supp}\mu$ is $1$-regular at $x$, and even that it is $C^1$-regular at $x$.  Indeed, if $(x_n)$ is a sequence of points of ${\rm supp}(\mu)$ that converges to $x$ and $v_n=(X_n, Y_n)\in \Cc_{x_n}({\rm supp}\mu)$   converges  to $v=(X,Y)$, we have for every $n$: 
 $$\| Y_n-\tilde s_-(x_n)X_n\|\leq 2\sqrt{\Delta s(x_n)}\sqrt{\Delta s(x_n)(X_n, X_n)}.$$
 As $G_-(x)=G_+(x)$, $\tilde s_-$ and $\Delta s$ are continuous at $x$. We deduce that 
 $\| Y-s_-(x)Y\|=0$ and then $(X, Y)\in G_-(x)$. We have then proved: 
\begin{cor}
If $\mu$ is an ergodic minimizing measure   all of whose Lyapunov exponents are zero, then, ${\rm supp}\mu$ is $C^1$ regular at $\mu$-almost every point. 
\end{cor}
This is exactly theorem \ref{th1}.
\section{Proof of the results of section \ref{sec3}}\label{sec4}

In this section, we   use the images of the physical verticals to obtain a control of the weak KAM solutions. More precisely, we can choose a graph in the image of a vertical, the graph of $da$ for a certain function $a$, and prove a certain inequality between $a$ and the considered weak KAM solution $u$. Then we deduce an inequality along some subset of the Ma\~n\'e set between the ``second derivatives'' of $a$ and $u$.  This gives a relationship between the Green bundles and the Bouligand's contingent cones to the pseudograph of any weak KAM solution along some subset of the Ma\~n\'e set .   
\subsection{Selection of some graphs in the images of the verticals}
\begin{nota}\begin{enumerate}
\item[$\bullet$] If $q\in M$, we denote the (physical) vertical $\pi^{-1}(\{ q\} )$ by $\Vc(q)\subset T^*M$.\\
\item[$\bullet$] If $t>0$, the function $A_t~: M\times M\rightarrow \R$ is defined by: 
$$A_t(q_0, q_1)=\inf_\gamma \int_0^tL(\gamma (s), \dot\gamma (s))ds=\min_\gamma \int_0^tL(\gamma (s), \dot\gamma (s))ds$$
where the infimum is taken on the set of $C^2$ curves  $\gamma~: [0, t]\rightarrow M$ such that $\gamma (0)=q_0$ and $\gamma (t)=q_1$. 
\item[$\bullet$] if $u~: M\rightarrow \R$ is a Lipschitz function, then by Rademacher's theorem, it is differentiable (Lebesgue) almost everywhere  and the   graph of its derivative is denoted by:
$$\Gc (du)= \{ (q, du(q)); u\ {\rm is} \ {\rm differentiable}\ {\rm at}\  q\} .$$
\end{enumerate}
\end{nota}

Tonelli's theorem asserts   that for every $t\not=0$, $\pi\circ\varphi_t(\Vc(q))=M$ (i.e. for every $q'\in M$ there exists a solution to the Euler-Lagrange equations $\gamma$ such that $\gamma (0)=q$ and $\gamma (t)=q'$); but in general $\varphi_t(\Vc(q))$ is not a graph. To select a  graph in $\varphi_t(\Vc(q))$,   we prove:
\begin{prop}
Let $H~: T^*M\rightarrow \R$ be a Tonelli Hamiltonian and $L~: TM\rightarrow \R$ be the associated Lagrangian. Then for every $t>0$ and every $q\in M$, the function $v_q^t=A_t(q,.)$ and $v_q^{-t}= A_t(.,q)$ are semi-concave, and satisfy:
$$\Gc(dv_q^t) \subset \varphi_t(\Vc(q))\ {\rm and}\  \Gc(-dv_q^{-t})\subset \varphi_{-t}(\Vc (q)).$$

\end{prop}

\demo
Because  $A_t$ is semi-concave, the  two functions $v^t_q $ and $v^{-t}_q$ are semi-concave and then Lipschitz.  By Rademacher's theorem  they are differentiable almost everywhere.\\
  Moreover, if $q_0$ is a point where $v^t_q$ is differentiable, then $v^t_q$ has exactly one super-differential at this point, there is  only one minimizing arc $\gamma$ joining $(0, q)$ to $(t, q_0)$, and we have: 
\begin{enumerate}
\item[$\bullet$] $dv_q^t(q_0)=\frac{\partial L}{\partial v}(\gamma (t), \dot\gamma (t))$;
\item[$\bullet$] $(\gamma(0),  \frac{\partial L}{\partial v}(\gamma (0), \dot\gamma (0)))=(q,  \frac{\partial L}{\partial v}(\gamma (0), \dot\gamma (0)))\in \Vc(q)$;
\item[$\bullet$] $\varphi_t \left( q,  \frac{\partial L}{\partial v}(\gamma (0), \dot\gamma (0))\right)=(\gamma (t), \frac{\partial L}{\partial v}(\gamma (t), \dot\gamma (t))) =(q_0, dv_q^t(q_0))$.
\end{enumerate}
Then we have proved that: $\varphi_t(\Vc(q))\supset \Gc (dv_q^t)$. Hence, we have selected a pseudograph in the image $\varphi_t(\Vc(q))$ of the vertical.\\
 In a very similar way, we may see that the  anti-pseudograph of the semi-convex function $-v_q^{-t}$ is a subset of $\varphi_{-t}(\Vc(q))$: $  \Gc (-dv_q^{-t})\subset \varphi_{-t}(\Vc(q))$.
 \enddemo
 
 \subsection{Local smoothness of some of these graphs  }
  \begin{nota}
 For every $x\in T^*M$, we denote  the {\em linear} vertical at $x$ by $V(x)$: $V(x)=\ker D\pi (x)=T_x\Vc(\pi(x))\subset T_x(T^*M)$.\\
The images of the linear vertical are denoted by: $G_t(x)=D\varphi_tV(\varphi_{-t}x)$.
 \end{nota}

 We recall that an orbit piece $(\varphi_t(x))_{t\in[a, b]}$ has no conjugate vectors if:
 $$\forall s\not=t\in[a,b], G_{t-s}(\varphi_tx)\cap V(\varphi_tx)=D\varphi_{t-s}(V(\varphi_s(x))\cap V(\varphi_t(x))=\{ 0\}.$$
  \begin{nota}
  Let us now fix a minimizing arc $\gamma~: [-t, 0]\rightarrow M$ such that:
  \begin{enumerate}
  \item[$\bullet$] there is only one minimizing arc between $(-t, \gamma(-t))$ and $(0, \gamma(0))$ (then it is  $\gamma$);
  \item[$\bullet$] the orbit piece $\left( \gamma(\tau), \frac{\partial L}{\partial v}(\gamma (\tau),\dot\gamma (\tau))\right)_{\tau\in[-t, 0]}$ has no conjugate vectors.
  \end{enumerate}
  Let us notice that when $(q, p)\in\Nc^*(H)$, then any piece of the curve  $(t\rightarrow \pi\circ\varphi_t(q,p))$ satisfies the previous hypotheses.\\
     We define a function  $a_t^+~: M\rightarrow \R$  by: $a_t^+(q)=v_{\gamma(-t)}^t(q)=A_t(\gamma(-t), q)$ (this function depends on $\gamma$). \\
     In a similar way, we can consider $x_0=(q_0, p_0)$ such that the orbit $(\varphi_s(x_0))_{s\in[0, t]}$ has no conjugate points and so that   there is only one minimizing arc $\gamma~: [0, t]\rightarrow M$ joining $q_0$ to $q_t$. We  define a  function   $ a_t^- ~: M\rightarrow \R$ by: 
  $  a_t^-=-v^{-t}_{q_{t}}(q)=-A_t(q, q_t)$. 
 \end{nota}
 
 \begin{prop}
 Let $\gamma~: [-t, 0]\rightarrow M$ (resp. $\gamma~: [0, t]\rightarrow M$) be a minimizing arc such that:
 \begin{enumerate}
 \item[$\bullet$] $\gamma$ is the only minimizing arc joining its two ends;
 \item[$\bullet$] the orbit piece $(\gamma, \frac{\partial L}{\partial v}(\gamma, \dot\gamma))$ has no conjugate vectors.
 \end{enumerate}
 Then there exists a neighborhood $V_0$ of $q_0=\gamma (0)$ in $M$ such that $a^+_{t|V_0}$ (resp. $a_{t|V_0}^-$) is as regular as $H$ is (then at least $C^2$).
 \end{prop}
 
 \demo
 We have seen   that:
$\Gc (da_t^+)\subset \varphi_t(V(q_{-t}))$. Let us now prove that   $a_t^+$ is smooth near $q_0$.

We use now   the so-called ``a priori compactness lemma'' (see \cite{Fa1}) that says to us that there exists a constant $K_t=K>0$ such that  the velocities $(\dot\gamma (s))_{s\in[0, t]}$ of any minimizing arc between any points $q\in M$ and $q'\in M$ are bounded by $K$; hence if we denote the set of the minimizing arcs that are parametrized by $[0, t]$ by $\Kc$, $\Kc$ is  a compact set for the $C^1$ topology because it is the image by the projection $\pi$ of a closed set of bounded orbits. Let us denote the set of $\gamma\in \Kc$ such that $\gamma (0)=q_{-t}$ by $\Kc_0$; then $\Kc_0$ is compact. Let us introduce another notation: 
$\Kc (q)=\{ \gamma\in \Kc_0; \gamma (t)=q\}$. Then $\Kc (q_0)=\{\gamma_0\}$ and hence, because $\Kc_0$ is closed,  for $q$ close enough to $q_0$, all the elements of $\Kc (q)$ are $C^1$ close to $\gamma_0$.

Moreover, $\varphi_t(\Vc(q_{-t}))$ is a sub-manifold of $M$ that contains $(q_0, \frac{\partial L}{\partial v}(q_0, \dot\gamma_0(0)))=(q_0, p_0)$. Its tangent space at $(q_0, p_0)$ is $G_t(q_0, p_0)$, which is transverse to the vertical  because $(q_s,p_s)_{s\in [-t, 0]}$ has no conjugate vectors. Hence, the manifold $\varphi_t(\Vc(q_{-t}))$ is, in a neighborhood $U_0$ of $(q_0, p_0)$, the graph of a $C^1$ section of $T^*M$ defined  on a neighborhood $V_0$ of $q_0$ in $M$.  Moreover, because this sub-manifold is Lagrangian (indeed, $\Vc (q_{-t})$ is Lagrangian and $\varphi_t$ is symplectic), it is the graph of $du_0$  where $u_0~: V_0\rightarrow \R$ is a $C^2$ function. 

Now, if $q$ is close enough to $q_0$, we know that all the elements $\gamma$ of $\Kc (q)$ are $C^1$ close to $\gamma_0$, and then that $(q, \frac{\partial L}{\partial v}(\gamma (t), \dot\gamma (t)))$ belongs to the neighborhood $U_0$ of  $(q_0, p_0)= (q_0, \frac{\partial L}{\partial v}(\gamma_0 (t), \dot\gamma_0 (t)))$ and to $\varphi_t(\Vc(q_{-t}))$. Because $\varphi_t(\Vc(q_{-t}))\cap U_0$ is a graph,   this element is unique: $\Kc (q)$ has only one element and $a_t^+$ is differentiable at $q$, with $da_t^+(q)= \frac{\partial L}{\partial v}(\gamma (t), \dot\gamma (t)))=du_0(q)  $ . We deduce that near $q_0$, on the set of differentiability of $a_t^+$, $da_t^+$ is equal to $du_0$; because $a_t^+$ and $u_0$ are Lipschitz on $V_0$ and their differentials are equal almost everywhere, we deduce that on $V_0$, $a_t^+-u_0$ is constant. Hence, on a neighborhood $V_0$ of $q_0$, $a_t^+$ is $C^2$. 

In a similar way, using the fact that $ \Gc(da_t^-)\subset \varphi_{-t}(V(q_t))$, we obtain that $a_t^-$ is $C^2$ near $q_0$.
\enddemo

\begin{remk}
If $x_0=(q_0, p_0)$ is a point of the Ma\~n\'e set, $(q_t,p_t)_{t\in\R}=(\varphi_t(q_0, p_0))_{t\in\R}$ has no conjugate vectors and for every $t<\tau$, there is only one minimizing arc $\gamma~:[t, \tau]\rightarrow M$ joining $q_t$ to $q_\tau$, hence for every $t>0$ the two functions $a_{x_0, t}^+$ and $a_{x_0, t}^-$ are smooth near $q_0$ (of course the neighborhood of $q_0$ where they are smooth depends on $t$).
\end{remk}

\subsection{Comparison between the weak K.A.M. solutions and the maps $a_t^+$ and $a_t^-$}
  
  \begin{lemma}
  We assume that $u_-$ is a negative  weak K.A.M. solution and that $u_+$ is a positive weak K.A.M. solution. Let $q_0\in M$ be a point of differentiability of $u_-$ (resp. $u_+$) and $a_t^+$ (resp. $a_t^-$) be the function built in the previous subsection for the arc $\gamma=(\pi\circ \varphi_s(q_0, du_-(q_0)))_{s\in [-t, 0]}$ (resp. $\gamma=(\pi\circ \varphi_s(q_0, du_+(q_0)))_{s\in [0, t]}$). Then, in a chart:  $u_-(q)-u_-(q_0)-du_-(q_0)(q-q_0)\leq a_t^+(q)-a_t^+(q_0)-da_t^+(q_0)(q-q_0)$ 
 (resp. 
 $a_t^-(q)-a_t^-(q_0)-da_t^-(q_0)(q-q_0)\leq u_+(q)-u_+(q_0)-du_+(q_0)(q-q_0)$).

  \end{lemma}
  \demo
  Let us consider $q_0$ in $M$ that is a point of differentiability of a weak K.A.M. solution $u_-$ and let us denote the point above $q_0$ on the pseudograph $\Gc (du_-)$ of $u_-$ by $x_0$: $x_0=(q_0, du_-(q_0))$. Then, for every $t>0$, because $T_tu_-=u_--ct$ is differentiable at $q_0$,   there is only one point $q\in M$ such that $u_-(q_0)=T_tu_-(q_0)+ct=u(q)+A_t(q, q_0)+ct$  and only one 
 minimizing arc $\gamma~: [-t, 0]\rightarrow M$ joining  $q$ to $q_0$.  We introduce the notation: $x_t=(q_t, p_t)=\varphi_t(x_0)$. Then: $T_tu_-(q_0)=u_-(q_{-t})+A(q_{-t}, q_0) $; moreover:  $T_tu_-(q)\leq u_-(q_{-t})+A(q_{-t}, q) =T_tu_-(q_0)+A(q_{-t}, q)-A(q_{-t}, q_0)$.  Finally: $u_-(q)-u_-(q_0)\leq a_t^+(q)-a_t^+(q_0)$. Because these two maps $a_t^+$ and $u_-$ are differentiable at $q_0$, they have the same differential at this point and  we obtain (in chart):
 $u_-(q)-u_-(q_0)-du_-(q_0)(q-q_0)\leq a_t^+(q)-a_t^+(q_0)-da_t^+(q_0)(q-q_0)$.
  
  Using the same argument for $u_+$, we obtain:
 if $q_0$ is a point of differentiability of $u_+$: \\
 $a_t^-(q)-a_t^-(q_0)-da_t^-(q_0)(q-q_0)\leq u_+(q)-u_+(q_0)-du_+(q_0)(q-q_0)$.
 
 \enddemo
 
 \bigskip

Now we would like to use these inequalities at different points $q_0$; we have to be careful, because $a_t^+$ and $a_t^-$ depend on the point $q_0$ we choose. That is why we change now our notation, replacing $a_t^+$ by $a_{q_0, t}^+$  if the considered point is $(q_0, du_-(q_0))$ and $a_t^-$ by $a_{q_0, t}^-$  if the considered point is $(q_0, du_+(q_0))$. 

\begin{prop}
We assume that $u_-$ is a negative  weak K.A.M. solution and that $u_+$ is a positive weak K.A.M. solution. Let $y\in\Ic(u_-, u_+)$ be a point, $(x_n)$ be a sequence of points of $M$ converging to $y$, and $(t_n)$ be a sequence of positive real numbers so that the two limits (written in charts) $\displaystyle{\lim_{n\rightarrow \infty} \frac{x_n-y}{t_n}=X}$ and $\displaystyle{Y=\lim_{n\rightarrow \infty} \frac{du_-(x_n)-du_-(y)}{t_n} }$ (resp. $\displaystyle{ \lim_{n\rightarrow \infty} \frac{du_+(x_n)-du_+(y)}{t_n} }$)  exist. Then we have:
$$\forall k\in\R^n, Y.k \leq \frac{1}{2}(d^2a_{y,t}^+(y)(k,k)+d^2a_{y,t}^+(y)(X,X)-d^2a_{y,t}^-(y)(X-k,X-k))$$
(resp:
$$\forall k\in\R^n,  \frac{1}{2}(d^2a_{y,t}^-(y)(k,k)+d^2a_{y,t}^-(y)(X,X)-d^2a_{y,t}^+(y)(k-X,k-X))\leq Y.k)$$
\end{prop}

\demo
We work in a chart, and we have, if $ y\in \Ic(u_-, u_+)$ and $x$ is a point of differentiability of $u_-$:
\begin{enumerate}
\item[$\bullet$] $u_-(x+h)-u_-(x)-du_-(x)h\leq a_{x,t}^+(x+h)-a_{x,t}^+(x)-da_{x,t}^+(x)h$;
\item[$\bullet$] $u_-(x)-u_-(y)-du_-(y)(x-y)\leq a_{y,t}^+(x)-a_{y,t}^+(y)-da_{y,t}^+(y)(x-y)$;
\item[$\bullet$] $a^-_{y,t}(x+h)-a^-_{y,t}(y)-da_{y,t}^-(y)(x+h-y)\leq u_+(x+h)-u_+(y)-du_+(y)(x+h-y)$.
\end{enumerate}
Hence, by adding these three inequalities and  using that $u_-(y)=u_+(y)$, $du_-(y)=du_+(y)$  and $u_+\leq u_-$:\\

\noindent$(du_-(y) -du_-(x))h   \leq a_{x,t}^+(x+h)-a_{x,t}^+(x)-da_{x,t}^+(x)h+a_{y,t}^+(x)-a_{y,t}^+(y)-da_{y,t}^+(y)(x-y) -a^-_{y,t}(x+h)+a^-_{y,t}(y)+da_{y,t}^-(y)(x+h-y)$.\\

We now   need to precise the regularity of the maps: $x\rightarrow da^-_{x,t}$  and $x\rightarrow da^+_{x,t}$. To do that, we prove a lemma.  We fix  a finite atlas of $M$ to write that $u_-$ is $K$-semi-concave and that $u_+$ is $K$-semi-convex. The proof is very similar to the one given by A.~Fathi in \cite{Fa1} to prove that the Aubry set is a Lipschitz graph.

\begin{lemma}\label{Aubrylipschitz}There exists a constant $K>0$  such that, for every $y\in\Ic(u_-, u_+)$ and every $x\in M$ where $u_-$ (resp. $u_+$) is differentiable, then $\| du_-(y)-du_-(x)\|\leq K\| y-x\|$ (resp. $\| du_+(y)-du_+(x)\|\leq K\| y-x\|$ ). In particular, $du_-$ and $du_+$ are continuous at every point of $\Ic(u_-, u_+)$.
\end{lemma}
\demo Because $u_+\leq u_-$, $u_-$ is semi-concave and $u_+$ is semi-convex, then  $u_-$ is   $K$ semi-convex at every point of $\Ic(u_-, u_+)$; hence: 
\begin{enumerate}
\item[$\bullet$] $u_-(x+h)-u_-(x)-du_-(x)h\leq K\| h\|^2$;
\item[$\bullet$] $u_-(x)-u_-(y)-du_-(y)(x-y)\leq K\|x-y\|^2$;
\item[$\bullet$] $-K\| x+h-y\|^2 \leq u_-(x+h)-u_-(y)-du_-(y)(x+h-y)$.
\end{enumerate}
Adding these three inequalities, we obtain: \\
 $(du_-(y)-du_-(x))h\leq K\| h\|^2+K\|x-y\|^2+K\| x+h-y\|^2$.\\
  We choose $h$ such that $\| h\|=\| x-y\|$: 
  $(du_-(y)-du_-(x))\frac{h}{\| h\|}\leq 6K\| x-y\|$ and then: $\| du_-(x)-du_-(y)\|\leq 6K\| x-y\|$. We have found a constant for $y$ close to $x$, this is enough to conclude because $\Ic (u_-, u_+)$ is compact and $du_-$ is bounded on $M$.
  \enddemo

  Let us now fix $y\in\Ic(u_-, u_+)$. For $x$ close to $y$ that is a point of differentiability of $u_-$, we have: 
 \begin{enumerate}
 \item[$\bullet$]$a_{x,t}^+(z)=A_t(\pi\circ\varphi_{-t}(x, du_-(x)), z)$; 
 \item[$\bullet$]$\{ (z, da_{x,t}^+(z))\}=\Vc (z)\cap  \varphi_t(\Vc_{\rm loc}(\pi\circ\varphi_{-t}(x, du_-(x))))$; \item[$\bullet$]
  ${\rm graph} (d^2a_{x,t}^+(z))  $ $=
   $ $T_{(z, dg_{x,t}^+(z))}D\varphi_t(\Vc (\pi\circ\varphi_{-t}(x, du_-(x))))=G_t(x, du_-(x))$ and then the previous intersection is transverse. 
   \end{enumerate}These three quantities depend  on $x$ and $z$; because $du_-$ is continuous at $y$, we have: for every $\varepsilon>0$, there exists $\delta>0$ such that, if $\| x-y\|<\delta$ and $z$ is in the chart near $y$:
  $\|d^2a^+_{x,t}(z)-d^2a_{y,t}^+(z)\|\leq \varepsilon$. \\
  Moreover, by Taylor-Lagrange inequality, we have:\\
 $\displaystyle{ \|a_{x,t}^+(x+h)-a_{x,t}^+(x)-da_{x,t}^+(x)h-\frac{1}{2}d^2a_{x,t}^+(x)(h,h)\|\leq\max_{z\in [x, x+h]}\| d^2a_{x,t}^+(z)-d^2a^+_{x,t}(x)\|\|h\|^2.}$
 Hence, if $x$ is close enough to $y$ and $h$ small enough~:\\
$\displaystyle{ \|a_{x,t}^+(x+h)-a_{x,t}^+(x)-da_{x,t}^+(x)h-\frac{1}{2}d^2a_{y,t}^+(x)(h,h)\|\leq \varepsilon \|h\|^2.} $\\
We have of course a similar result for $a_{x,t}^-$ and $x$ any differentiability point of $u_+$.\\
Let us now consider a sequence $(x_n)$ of points of differentiability of $u_-$ that converges to $y$ so that: $\forall n, x_n\not=y$,   a vector $k$ with fixed norm $\| k\|=\lambda>0$ and $(h_n)=(t_nk)$ where   $(t_n)$ is a sequence of positive numbers tending to $0$. We have proved that:

\noindent$(du_-(y) -du_-(x_n))h_n   \leq a_{x_n,t}^+(x_n+h_n)-a_{x_n,t}^+(x_n)-da_{x_n,t}^+(x_n)h_n+a_{y,t}^+(x_n)-a_{y,t}^+(y)-da_{y,t}(y)(x_n-y) -a^-_{y,t}(x_n+h_n)+a^-_{y,t}(y)+da_{y,t}^-(y)(x_n+h_n-y)$.\\
We assume that $\displaystyle{\lim_{n\rightarrow \infty} \frac{ x_n-y}{t_n}=X}$ and $\displaystyle{Y=\lim_{n\rightarrow \infty} \frac{du_-(x_n)-du_-(y)}{t_n} }$.

We divide by $t_n^2$ the previous inequality and take the limit when $n$ tends to $+\infty$ and we obtain: 
$$-Y.k \leq \frac{1}{2}(d^2a_{y,t}^+(y)(k,k)+d^2a_{y,t}^+(y)(X,X)-d^2a_{y,t}^-(y)(X+k,X+k))$$
changing $k$ into $-k$, this gives the wanted result. 
In a similar way we obtain for $u_+$:
$$\forall k\in\R^n,  \frac{1}{2}(d^2a_{y,t}^-(y)(k,k)+d^2a_{y,t}^-(y)(X,X)-d^2a_{y,t}^+(y)(k-X,k-X))\leq Y.k $$\enddemo

\subsection{Links between the Green bundles and the weak K.A.M. solutions}
\begin{nota} Near every point $q\in M$, we choose    some coordinates $(q_1, \dots , q_n)$ of $M$ and associate to them their dual coordinates $(p_1, \dots , p_n)$ such that $(q_1, \dots, q_n, p_1, \dots , p_n)$  are symplectic coordinates on $T^*M$. Then we can associate to these coordinates their infinitesimal coordinates $(Q_1, \dots, Q_n, P_1, \dots , P_n)$. \\
Then  any  Lagrangian subspace $G$ of $T_x(T^*M)$ that is transverse to the vertical is the graph of a linear map  whose matrix $s$ in the coordinates $(Q_1, \dots, Q_n, P_1, \dots , P_n)$ is symmetric. We can then associate to $G$ the unique quadratic form $Q$ whose matrix (as a quadratic form) in coordinates $(Q_1, \dots , Q_n)$ is $s$.

For example, if $q\in M$ is a point of  differentiability of $u_-$ (resp. $u_+$ then  the Green bundle $G_+(q, du_-(q))$ (resp. $G_-(q, du_+(q))$) is well  defined and transverse to the vertical. We denote  by $Q_-$ (resp. $Q_+$) its associated quadratic form and by $s_-$ (resp. $s_+$) its matrix.  \end{nota}
Let us recall that if $x\in A\subset T^*M$, $\Cc_xA$ designates the contingent cone to $A$ at $x$, that was defined in the introduction.
\begin{prop}
We assume that $(u_-, u_+) $ is a pair of conjugate weak K.A.M. solutions. Let $y\in\Ic(u_-, u_+)$ be a point  and $(X,Y)\in \Cc_{(y, du_-(y))}\Gc(du_-)$. Then we have:
$$\forall k\in\R^n, Y.k   \leq \frac{1}{2}(Q_+(k,k)+Q_+(X,X)-Q_-(X-k,X-k)) $$
and if $(X,Y)\in \Cc_{(y, du_+)}\Gc (du_+)$:
$$\forall k\in\R^n, \frac{1}{2}(Q_-(k,k)+Q_-(X,X)-Q_+(X-k,X-k))\leq Y.k .$$

\end{prop}
\demo
We know   that $\displaystyle{G_+(q,p)=\lim_{t\rightarrow +\infty}G_t(q,p)}$ (resp. $\displaystyle{G_-(q,p)=\lim_{t\rightarrow -\infty}G_t(q,p)}$). Hence, if $q$ is a point of differentiability of $u_-$, we have:  $\displaystyle{Q_+(q, du_-(q))=\lim_{t\rightarrow +\infty} d^2g_{q,t}^+(q)}$ and if $q$ is a point of differentiability of $u_+$: $\displaystyle{Q_-(q, du_+(q))=\lim_{t\rightarrow +\infty} d^2g_{q,t}^-(q)}$. If we use the inequalities that we proved in the previous section,  we obtain: 
$$\forall k\in\R^n,   Y.k   \leq \frac{1}{2}(Q_+(X,X)+Q_+(k,k)-Q_-(X-k,X-k)).$$

Let us now look for the contingent cone to the pseudograph $\Gc (du_-)$ at $(y, du_-(y))\in \Ic(u_-, u_+)$. Working in a chart, we assume that $(X, Y)\in \Cc_{y, du_-(y))}\Gc (du_-)$ is not the null vector. Hence, there exists a sequence $(t_n)$ of positive numbers that converges to $0^+$ and a sequence $(x_n)$ of points of differentiability of $u_-$ that converges to $y$ so that:
$$(X, Y)=\lim_{n\rightarrow \infty} \frac{1}{t_n}(x_n-y, du_-(x_n)-du_-(y)).$$
 This corresponds exactly to the limit that we computed in the previous subsection. Hence, we  proved:
\medskip

{\em If $y\in \Ic(u_-, u_+)$, if $(X, Y)$ is a vector of the contingent cone to $\Gc (du_-)$ at $(y, du_-(y))$, then:

$$\forall k\in \R^n,  Y.k\leq \frac{1}{2}(Q_+(k,k)+Q_+(X,X)-Q_-(X-k, X-k)).$$}
In a similar way, we obtain: 
\medskip

{\em If $y\in \Ic(u_-, u_+)$, if $(X, Y)$ is a vector of the contingent cone to $\Gc (du_+)$ at $(y, du_+(y))$, then:

$$\forall k\in \R^n,  \frac{1}{2}(Q_-(k,k)+Q_-(X,X)-Q_+(X-k, X-k)\leq Y.k .$$}
\enddemo

\subsection{Proof of theorem \ref{greenkam}}

Let $(u_-, u_+)$ be a pair of conjugate weak KAM solutions and let $q$ belong to $\Ic(u_-, u_+)$. We want to prove that: $\forall (X, Y)\in  \Cc_{(q, du_-(q))}\Gc (du_-),  $ 
$$\| Y-\tilde s_-(q, du_-(q))X\|\leq 2\sqrt{\|\Delta s(q, du_-(q))\|} .\sqrt{\Delta s(q, du_-(q))(X,X)}$$
 $$\leq 2\Lambda (\Delta s(q, du_-(q))). \| p_{\Delta s(q, du_-(q))}(X)\|$$

 We denote the quadratic form associated with $G_-$ (resp. $G_+$) by $Q_-$ (resp. $Q_+$). Then the quadratic form associated with $\tilde G_-$ (resp. $\tilde G_+$) is $\tilde Q_-=2Q_--Q_+$ (resp. $\tilde Q_+=2Q_+-Q_-$). Let $(X, Y)\in\Cc_{(q, du_-(q))}\Gc (du_-)$ be a vector of the contingent cone. We have proved that: 
$$\forall k\in \R^n,  Y.k\leq \frac{1}{2}(Q_+(k,k)+Q_+(X,X)-Q_-(X-k, X-k)).$$
Then we write: $Y={}^tQ_+X+\Delta Y$ and $\Delta Q=Q_+-Q_-$. The previous inequality can be rewritten as follows:
$$\forall k\in \R^n,  \Delta Y.k\leq \frac{1}{2}\Delta Q(X-k, X-k)\quad (*).$$

We have the following splitting: $\R^n=\ker {}^t \Delta Q\oplus \Im {}^t\Delta Q$ and $\Delta Y=Y_1+Y_2$ with $Y_1\in \ker {}^t\Delta Q$ and $Y_2\in \Im{}^t\Delta Q$. We deduce from $(*)$:
$$\forall k\in\ker{}^t\Delta Q, Y_1. k\leq  \frac{1}{2}\Delta Q(X, X). $$
This implies: $Y_1=\vec 0$. Hence $\Delta Y=Y_2\in \Im {}^t\Delta Q$ and there exists a unique $y\in -2X+\Im {}^t\Delta Q$ such that $\Delta Y={}^t \Delta Q y$. Then $(*)$ becomes:
$$\forall k\in\R^n, \Delta Q(y, k)\leq \frac{1}{2} \Delta Q(X-k, X-k)$$
i.e:
$$\forall k\in\R^n, \Delta Q(X+\frac{y}{2}, X+\frac{y}{2})-\Delta Q(X,X)\leq \Delta Q(X-k+\frac{y}{2}, X-k+\frac{y}{2})$$
As $G_-\leq G_+$, the quadratic form $\Delta Q$ is positive semi-definite. Hence the previous inequality is equivalent to:
$$\Delta Q(2X+y, 2X+y)\leq 4 \Delta Q(X,X).$$
 Let us write $y=-2X+\Delta y$. We have $\Delta y\in \Im {}^t\Delta Q$. Then 
$Y=  {}^t\tilde Q_-X+{}^t\Delta Q\Delta y$ and $\Delta Q(\Delta y, \Delta y)\leq 4 \Delta Q(X,X)$. \\
Then we can write: $\Delta y=2\sqrt{Q(X,X)}  u$ with $\Delta Q(u, u)\leq 1$ and:
$$\| {}^t\Delta Q\Delta y\|^2=4\Delta Q(X, X)({}^t\Delta Q u).({}^t\Delta Q u)$$
with:
$$({}^t\Delta Q u).({}^t\Delta Q u)\leq \sup \frac{{}^t\Delta Q v .{}^t\Delta Q v }{\Delta Q(v,v)}= \| \Delta Q\|.$$

We then obtain: $$\| Y-{}^t\tilde Q_-X\| \leq 2\sqrt{\| \Delta Q\|}.\sqrt{\Delta Q(X, X)}.$$
If we denote by $\Lambda(\Delta Q)$ the biggest eigenvalue of $\Delta Q$ and by $p_{\Delta Q}$ the orthogonal projection on the  image of ${}^t\Delta Q$, we deduce:
$$\| Y-{}^t\tilde Q_-X\| \leq 2\Lambda (\Delta Q) \| p_{\Delta Q}(X)\|Ê.$$

\enddemo

\end{document}